# Rejoinder: Classifier Technology and the Illusion of Progress

**David J. Hand**

I would like to thank the discussants for some very stimulating comments. Being only human, I am naturally pleased when others produce evidence or arguments in support of my contentions, but being a scientist, I am also pleased when others produce evidence or arguments against my proposals (although I may have to take a deep breath first), since this represents the scientific process in action.

I should first make one thing clear: I agree with Professor Friedman that substantial advances have been made in recent years. Indeed, in my paper I remarked that "developments such as the bootstrap and other resampling approaches ... have led to significant advances in classification and other statistical models." However, what I question is whether the advances, when taken in the context of real practical problems, are as great as is often claimed—the recognition of the limitations of the new methods to which Professor Friedman refers.

Professor Friedman agrees with my three points that the improvements of newer methods over older ones are less than those of the older ones over still older ones, that the evidence favoring the superiority of new methods is often suspect and that the new methods fail to tackle important problems. I draw the conclusion from these points that progress is not as great as is imagined. Professor Friedman draws the conclusion that low lying fruit is easier to gather, that initial validation of new methods should be more rigorous and that much work remains to be done. Perhaps, then, we are really broadly in agreement—only perhaps I am describing a half empty glass (the new classification tools are not as wonderful as they are claimed), while Professor Friedman is describing a half full glass (some classification tools represent advances over the older ones).

I admit that I did criticize error rate as a performance measure and then used it in the examples. Since most performance comparisons of classifiers use error rate, this seemed justifiable, and I believe that my conclusions will generalize to other performance measures. For example, I agree that in some two-class problems it is the rank order of the estimated class 1 membership probabilities which matters and that modern methods may well be able to estimate this more accurately than older methods. However, surely my points about population drift, class definition uncertainty and so on still apply and, of course, my point that people often use one criterion to fit a model and another to evaluate it applies even more strongly.

In fact, this point about people using different criteria manifests itself at a higher level when Professor Friedman and I examine my Table 1. I see the proportion of reduction of error rate achieved by the best method which can be achieved by discriminant analysis, whereas Professor Friedman sees the ratio of the error rates. I see a large initial improvement so that subsequent improvements are relatively small; he sees a large reduction in the proportion remaining. Back to the half full/half empty glasses again. We are both right, of course, although perhaps the different perspectives are valuable for different uses. For example, I agree with Professor Friedman's example of the zip code classifier—and here the ratio of error rates might be a sensible measure—but (I would imagine) this is a problem in which the distributions are fairly static. In other problems, the distributions will change rapidly and I can imagine many contexts when I would not want to place too much trust in a reduction of error rate by a factor even as large as 10, if it corresponded to a change from a starting point as small as 0.001 to an even smaller one of 0.0001. A slight shift in the shapes of the distributions might induce sufficiently large changes in error rate so as to make this change irrelevant.

My regression example in Section 2.1 was merely intended as an additional illustration of the fact that the sequential nature of modeling means that typically later improvements are smaller than early ones. I am suggesting that the first, relatively crude, models will generally yield greater marginal improvements in predictive power than the later models. This is the low hanging fruit phenomenon—athough, as noted below and as Professor Stine illustrates, there are exceptions.

I am glad Professor Friedman agrees so strongly with Section 5 of the paper, on the difficulties of obtaining generally valid empirical comparisons. I think this is one of the most important parts of the paper. Professor Friedman's suggestion that the top performers in comparison studies should be ignored and attention should be focussed on the relative rankings of the others is very valuable. I also recommend looking at those methods which generally perform well, even if they seldom perform best, since they will have some sort of robustness. I think these sorts of issues, which represent aspects of the *art* of statistics, are fundamental to good statistical practice. They are the sorts of things which are not taught in standard statistics texts.

Professor Friedman comments that certain methods (he uses ensemble methods and support vector machines as examples) "offer substantial advantages over the earlier methods in enough situations to be regarded as major advances." I agree that such methods do represent significant theoretical and practical advances. My point is the milder one that "the practical impact of the developments has been inflated; that although progress has been made, it may well not be as great as has been suggested." Again referring to population drift as an example, a better fit to data drawn from a given distribution is not so wonderful if the distribution has changed. In fact, of course, it is likely that Professor Friedman and I have slightly different experience in terms of application domains. He cites "scientific and engineering applications" and I cite examples such as credit scoring and fraud detection: he draws attention to the differences between domains toward the end of his contribution; it is possible that population drift is more apparent in the latter than the former.

I entirely endorse Professor Friedman's comment that "obtaining high quality representative training data is generally more important to success than choice of a particular classifier." We are agreed on this, but in part my paper aims to point out that obtaining "representative training data" may be harder than is often imagined. Incidentally, I often go one step further and suggest that the best way to dramatically improve classifier performance is to add suitably chosen extra discriminating variables—that this is likely to exceed the performance improvement attained by juggling with classification rules, but, of course, this does depend on the specifics of the application.

Professor Friedman points out that almost all modern procedures incorporate a regularization parameter that controls the goodness of fit to the training data, and that one way to overcome problems such as population drift or uncertainty in the class definitions is to regularize more heavily than one would if such problems were not suspected. I agree, and I also agree that there is no reason to suppose that the arbitrary amount of extra regularization implied by simpler older methods is the right amount. Indeed, of course one can always find examples where it is not, such as the large $d$ small $n$ cases of bioinformatics. However, if one is unable to get a handle on the amount of regularization which is needed, then there is no reason to suppose that the more heavily regularized modern method will be any better than the implicitly regularized older method.

Professor Friedman provides a useful discussion of tools for handling errors in class labels. These are fine if one suspects that one has such errors. However, I was concerned with the question of robustness to such errors if one is using a more standard method, unaware of the possibility.

I am sorry to have disappointed Professor Stine by not giving "a rich portfolio of examples that demonstrate the failures of complex models." To some extent I am caught in a Catch-22 situation here. For example, had I demonstrated the superiority of a simple linear classifier over a complex support vector machine in a real example involving dramatic population drift, then an obvious response would have been to build a more elaborate dynamic classifier or apply a modern model with heavier than standard regularization, as suggested by Professor Friedman. For this particular situation, the "even more elaborate model" would then win—and this will always be the case for any particular example. However, across examples, when one does not have specific reasons to expect such departures from the classic "fixed underlying distributions, precise class definitions" and so on of the standard problem, then one will not use a tool specially matched to the problem, so there is a risk that one will miss important features of the problem. Perhaps all I am really saying is that every problem has unique features, and that ideally one would carefully model and allow for those features, but if one is unaware of them (implicit in the use of standard tools), then simple is better.

My reason for using the idealized example of equally correlated predictors in Section 2.1 was merely to make the mathematics particularly transparent. Indeed, I pointed out that in real applications, the phenomenon I demonstrated was likely to be even more pronounced. However, I take Professor Stine's point that artificial examples can be used to support any argument (the intertwined spirals example being a case in point!), but, in spite of the ingenuity of his superadditive growth example, I believe that empirical evidence shows that decreasing marginal improvement as extra terms are added to a model is the norm.

I am not arguing that there are *no* contexts in which a small improvement in performance is valuable. Professor Stine's example of data compression is a nice one. Another, of course, would be a small improvement in classification accuracy in a medical screening context—correctly diagnosing people in time to be treated, for example. My argument relates this apparent small improvement to other sources of uncertainty in the problem. If the distributions of characteristics of people with the disease differ from the distributions used to construct the classification rule, then the apparent improvement may be illusory. Statistical significance does not affect this argument. If the distributions are not the right ones, it does not matter how statistically significant the apparent improvements are.

Like Professor Friedman, Professor Stine takes me to task for criticizing the use of inappropriate performance criteria (which we agree is wrong) but then using error rates in my example in Table 1. I agree, of course, and in an ideal world I would have used performance criteria better matched to the particular problems and objectives. To do this I would have had to use my own examples, for which I knew the relevant performance criteria, and then compared linear discriminant analysis with the best performance I could achieve using neural networks, support vector machines, random forests and the whole panoply of other methods. However, if I then tried to argue, as I did in the paper, that these sophisticated tools were not that much better than linear discriminant analysis, I would immediately be vulnerable to the criticism that this was simply because I was not very adept at using the other methods. I thought it would be more compelling to use the results of other, expert, analysts. This meant I was forced to use error rate in my comparisons, simply because this is the most widely used criterion.

Professor Stine's comment about the difficulty of extracting the full story from commercial clients, so that one is confident that one is answering the right question, struck a chord. Even worse, all too often the client is *incapable* of formulating a precise question. This is not intended as a criticism: often the intrinsic uncertainties of the world (especially the commercial world) make precise formulation impossible. This, of course, was one of the issues which stimulated my writing of the paper.

Professor Stine's example of population drift in a personnel selection problem is very nice. It involves the key issue of drift due to natural background changes (the economy), but also, presumably, the employees on which the model was built were not a random selection from previous applicants, but had been chosen because someone thought they were likely to be successful employees. This means, of course, that the classifier would have been modeling inappropriate distributions, unless some effort was made to represent this prior selection process. This is the same problem as that in the example of drivers of white cars which Professor Stine cites, although to a less extreme extent. I suspect that Professor Stine is right when he doubts that any model would have been very successful on this problem. Personnel selection problems are notoriously difficult. My point is merely that there are aspects of this problem which are not considered in the classical supervised classification paradigm, which consists of trying to model underlying distributions from a sample of data drawn from those distributions.

Toward the end of his contribution, Professor Stine asks for my suggestions on how to decide whether it is useful to look for extra structure. I think one should always look for this, but there are different kinds of structure. There is the structure represented by shape of the distribution from which the design data were drawn, and there is structure in the overall problem (e.g., population drift). I am suggesting that we are now pretty good at modeling the former, but that often the extra features of the distributions that our clever modern methods pick up are relatively unimportant compared with the potential impact of taking into account the latter kind of structure. So my answer to Professor Stine's question about my approach to deciding when additions to a simple model are worthwhile is that I think it is a matter of priorities. It is one thing to be able to add another hidden node to a neural network and hence reduce the misclassification rate (on those distributions) by 0.5%. It is another (and often a more useful) thing to be able to say that one is really interested in cost weighted error rate and is uncertain about the costs, so that the Gini coefficient is a more appropriate measure of performance, or that one believes the design data do not properly represent the distributions of new cases and so on.

As far as population drift is concerned, I think Professor Stine's final paragraph hits the nail on the head: statisticians now have a powerful armory of methods to tackle this, but how often does one see them integrated into the design of a classification rule?

It is in this vein that Professor Holte rightly points out that there are methods for dealing with some of the factors that I identify as being unknown at the time of classifier design or subject to change after that time. In fact, I would be surprised if methods do not exist for *all* such factors: the Kalman filters Professor Stine refers to for population drift, Heckman models for sample selectivity, the cost curves of Professor Holte and the weighted Gini coefficients of Adams and Hand for unknown relative misclassification costs, for example. In addition, if tools for coping with a particular kind of uncertainty in the problem indeed do not exist, then it is, as Professor Holte says, a challenge for future research. Even if such tools exist, how often are they applied? Once again I wonder if, perhaps, it is just that it is easier to refine an existing form of classification model (the extra nodes of the neural network, the more sophisticated metric in nearest neighbor methods, ...) than to model the sample distortion or adopt a more complicated performance criterion. Perhaps many of us academic researchers are still guilty of focusing too much on Tukey's exact answer to the wrong problem. I hope I may be forgiven for making that comment, since I, too, am an academic researcher and I, too, know the pleasure of developing a classification tool which appears to have a slight edge over its competitors.

Professor Gayler's comments were interesting, not least because they were from precisely the perspective which had stimulated many of my observations—the "nonclassical" problems which arise when applying supervised classifiers in the context of modeling human behavior, specifically credit scoring.

Professor Gayler points out the great financial gains which would result from a small increase in predictive accuracy in this application domain, so that one might have expected a premium to be placed on such performance, making the fact that relatively simple old-fashioned approaches are still used rather surprising. He also points out that the new methods are regularly investigated by the credit scoring community, but rarely make the transition to everyday practice, suggesting that the simpler older methods have some kind of advantage. Professor Gayler and I agree that this advantage arises from the kinds of issues described in my paper.

Professor Gayler mentions yet other kinds of complications. For example, he refers to account management changes (which will occur after the accept/reject classification has been made). This is a special case of a more general class of problem. Often we want to predict into what class an object (often a person) will fall if we take some action. However, if our prediction suggests that they will fall into some undesirable class, then we take some other action. This, of course, invalidates the prediction. It is a generalization of the reject inference problem and leads to particular sample selectivity issues.

Professor Gayler is right to point out that in many problems the value of the threshold (to be compared with the estimated probability of belonging to class 1, e.g.) above which objects are assigned to class 1 depends on operational decisions, and these will be determined by all sorts of external factors.

I was particularly struck by Professor Gayler's observation that "in the limit (and the hands of a skilled modeler), every modeling technique should end up in agreement because they are all approximating the same data." I am reminded of Hoadley's ping-pong theorem, which presumably represents alternate steps toward this limit! I was also taken by his suggestion that it might be more useful "to look at the effort required of the modeler to achieve a given goodness of fit and other properties of the models that are of operational relevance to the lender." I endorse this. Of what good (at least in the credit scoring context) is a tool so highly sophisticated that it can be used effectively only after years of practice and experience? Operational relevance is a key factor.

In fact, my comment about "intertwined spirals or checkerboard patterns" refers to more than problems which can be modeled only as interactions between the variables. I meant it also to refer to those problems which have an extremely complicated (or perhaps contrived) decision surface. Such problems appear to be extremely rare in the real world, so demonstration of the power of new methods by showing that they can tackle such problems is rarely relevant to real problems. I conjecture that such problems are rare because in real problems the predictor variables will generally have been chosen because they are thought to have some discriminatory power, and predicting that the classes would be separated in such a complex way by a combination of variables would be an extraordinary intellectual feat. It is much easier to identify variables on which the members of one class have a tendency toward higher values than the members of the other class.

I like Professor Gayler's observation that it would be undesirable (in a credit scoring context) for a small change in decisions made when modeling to lead to a large change in the models. This is true and is a nice example of the pressures that favor simple modeling strategies. In such an environment, the organizations need to be confident of their modeling strategy and that it will be reliable in the hands of other, perhaps less experienced staff. This is a phenomenon similar to, but at a level different from, the flat maximum effect. There the users of the models want to be confident that slight changes in the model (and indeed, the modeling conditions) will not lead to sudden dramatic deterioration in performance: as Gayler says, the flat maximum effect is a great advantage in credit scoring.

I still have a suspicion that there is too much emphasis on trying to squeeze the last drops of performance out of classifiers matched to a particular data set when these distributions might not be the right ones, when the performance criterion being used is inappropriate, when the class definitions might be incorrect or subject to change and so on, with all the mismatches illustrated in the paper and others. Instead, I believe that more effort should be spent on trying to identify and model aspects of the problem which deviate from the classical supervised classification paradigm, and which may have a substantial impact on performance. For example, if you suspect the populations will change (perhaps not in Professor Friedman's scientific and engineering problems, but certainly in the personnel and social applications of Professor Stine, Professor Gayler and myself), then either model this or regularize more heavily to allow for it; if you suspect that the sample has not been randomly drawn but has been purposively selected (as in Professor Stine's employee selection example), use a model which adjusts for the hypothesized selectivity or more heavily regularize to avoid overfitting a suspected inaccurate distribution; if you know you are concerned with maximizing profit, then use profit as a performance criterion, and not misclassification rate or likelihood, or else regularize more heavily to allow for the fact that there is a mismatch between the criterion being used and the one of real interest; and so on.

I am extremely grateful to the discussants for their thoughtful comments on the paper. It is apparent that they spent a considerable amount of time and effort carefully considering my points, and marshalling coherent and instructive responses. Their comments covered a wide range of issues and approached things from different perspectives. It is very clear that, whatever the merits of the paper itself, the discussion contributions have substantial intrinsic value, and I have certainly learnt a great deal from them.